\documentclass[letterpaper,11pt]{article}
\usepackage{amsmath, amsthm, amssymb}
\usepackage{bridges}
\usepackage{graphicx}
\usepackage[colorlinks=true, urlcolor=blue, citecolor=black, linkcolor=black]{hyperref}
\usepackage{subcaption}
\usepackage{float}

\usepackage{xcolor}

\title{Illustrating Hyperbolic Surfaces with Mesh Embeddings}

\author{Fabian Lander\textsuperscript{1}, Erik L\"offelholtz\textsuperscript{1}, Diaaeldin Taha\textsuperscript{1}, Steve Trettel\textsuperscript{2}, and Anna Wienhard\textsuperscript{1}
\vspace{10pt}\\
\textsuperscript{1}Max Planck Institute for Mathematics in the Sciences; \href{mailto:lander@mis.mpg.de}{\textcolor{black}{lander@mis.mpg.de}}\\
\textsuperscript{1}Max Planck Institute for Mathematics in the Sciences; \href{mailto:loffelholz@mis.mpg.de}{\textcolor{black}{loffelholz@mis.mpg.de}}\\
\textsuperscript{1}Max Planck Institute for Mathematics in the Sciences; \href{mailto:taha@mis.mpg.de}{\textcolor{black}{taha@mis.mpg.de}}\\
\textsuperscript{2}University of San Francisco; \href{mailto:strettel@usfca.edu}{\textcolor{black}{strettel@usfca.edu}}\\
\textsuperscript{1}Max Planck Institute for Mathematics in the Sciences; \href{mailto:wienhard@mis.mpg.de}{\textcolor{black}{wienhard@mis.mpg.de}}
}

\date{}

\begin{document}

\maketitle

\thispagestyle{empty}

\begin{abstract}

Hyperbolic geometry exhibits geometric phenomena, such as fast area growth, that are difficult to visualize faithfully in Euclidean space, and which standard models like the Poincar\'{e} disk can obscure. To bring hyperbolic geometry to life, we embed hyperbolic surfaces in Euclidean space by discretizing the surfaces into meshes, and minimizing a distortion energy so that the edge lengths in the embeddings match those in the hyperbolic plane. The resulting surfaces buckle and ruffle to accommodate the extra area, making visible what flat models hide. We present exemplary illustrations, such as embedded disks, equidistant strips, diverging geodesics, and also artistic organic-like renders.  We discuss our use of these models, as renders and 3D prints, in research talks, public engagement, outreach, and education.

\end{abstract}

\begin{figure}[h!tbp]
	\centering
	\includegraphics[width=\textwidth]{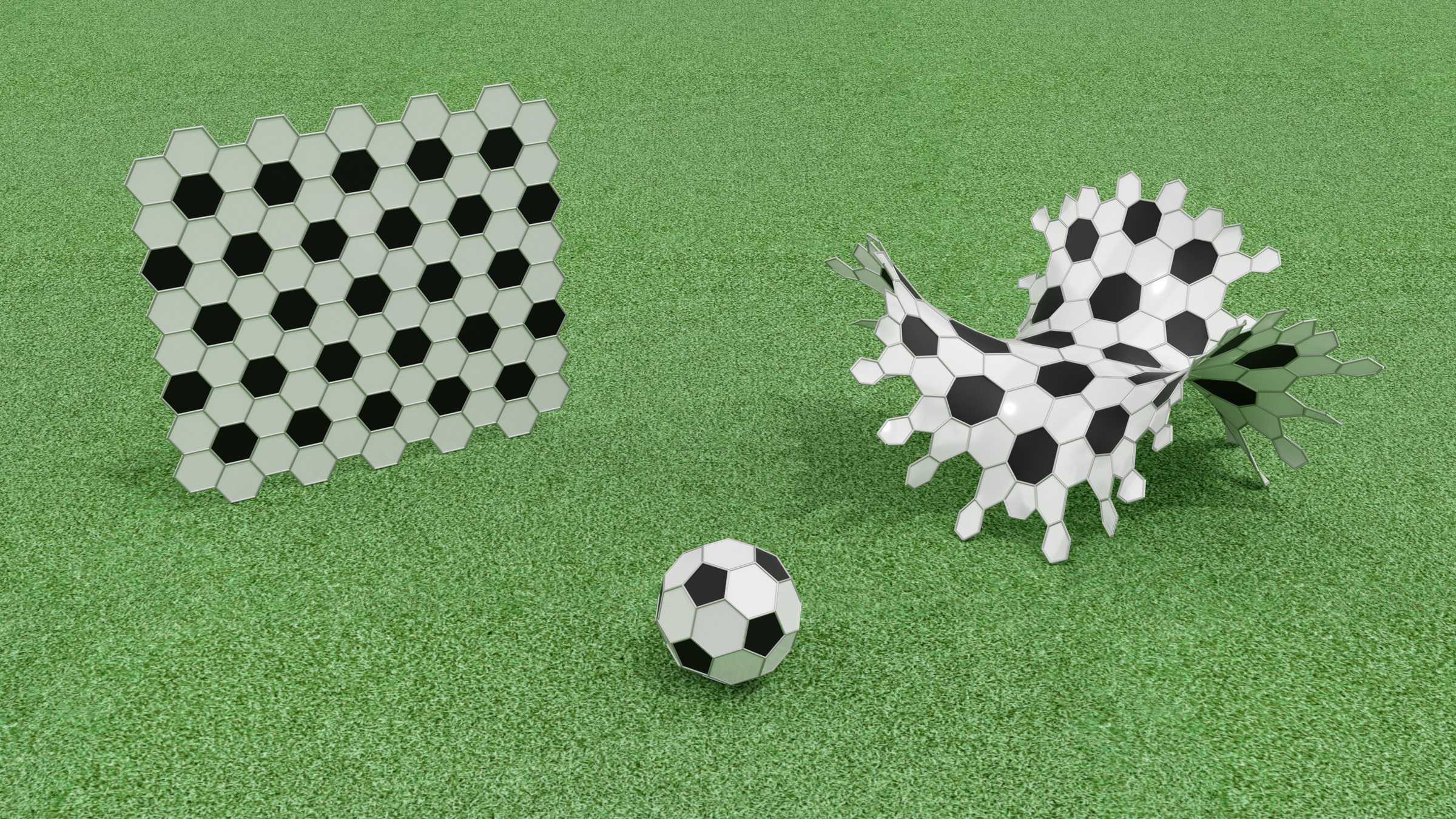}
	\caption{Three geometries, tiled in the pattern of a soccer ball. From left to right: Euclidean, spherical, and hyperbolic.}
	\label{fig:soccer}
\end{figure}

\section*{Introduction}

There are three geometries in which every point looks the same: the sphere, the Euclidean plane, and the hyperbolic plane. Figure~\ref{fig:soccer} shows all three, tiled in the pattern of a soccer ball. On the sphere, pentagons and hexagons fit together and close up. In the Euclidean plane, hexagons alone tile perfectly, extending forever. In the hyperbolic plane, this pattern continues, from pentagons to hexagons to heptagons, each step adding more area than in the Euclidean plane. Where Euclidean tilings grow flat, hyperbolic ones must buckle.

This extra area grows exponentially with distance, making hyperbolic geometry difficult to visualize. Unlike the sphere, which embeds naturally in $\mathbb{R}^3$, there are no known simple formulas for embedding regions of the hyperbolic plane into Euclidean space. The hyperbolic tiling in Figure~\ref{fig:soccer} has no closed-form embedding. Computing embeddings numerically for such illustrations is the subject of this paper.

\begin{figure}[h!tbp]
\centering
\begin{minipage}[b]{0.32\textwidth}
	\includegraphics[width=\textwidth]{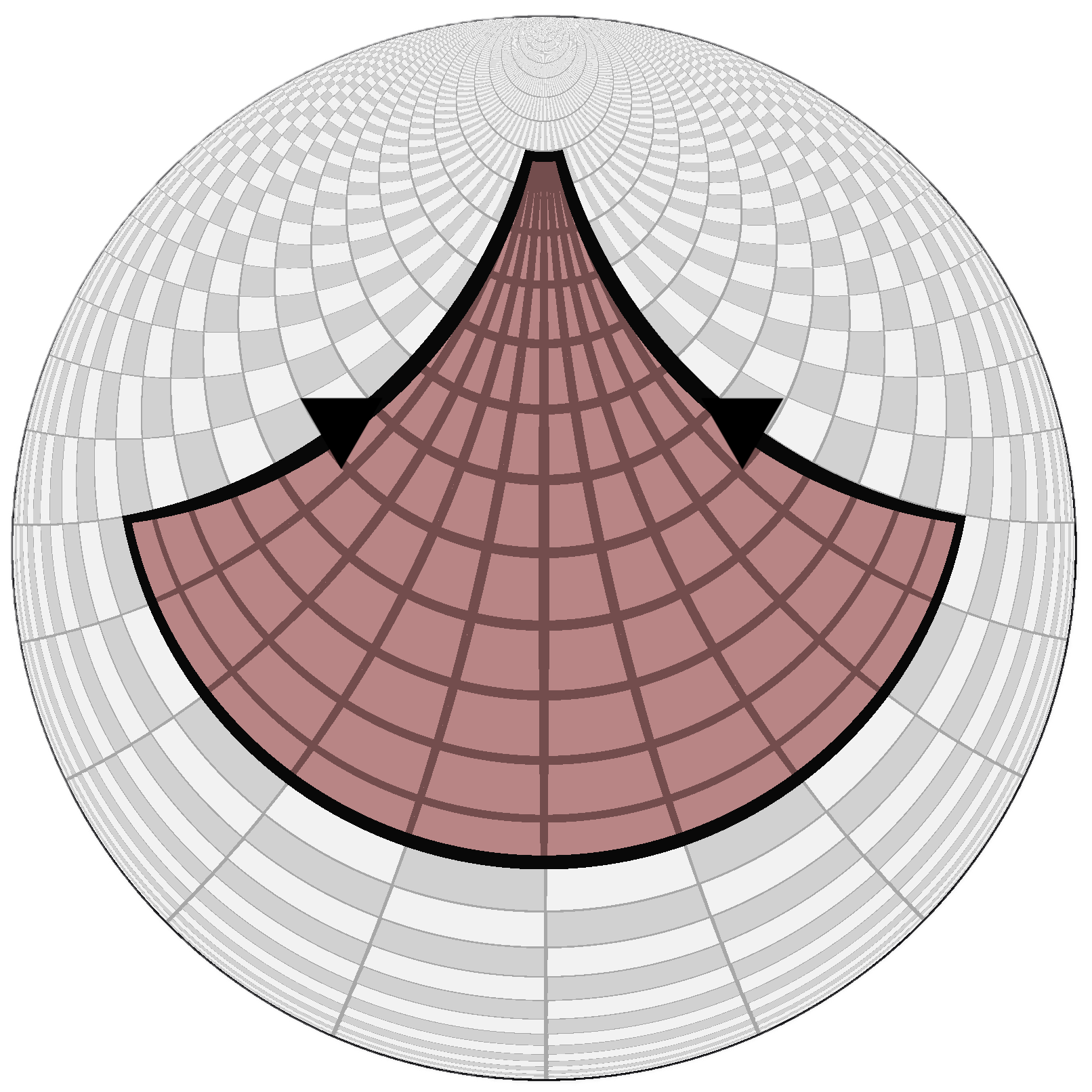}
        	\subcaption{}
        	\label{fig:1a}
\end{minipage}
~
\begin{minipage}[b]{0.32\textwidth}
	\includegraphics[width=\textwidth]{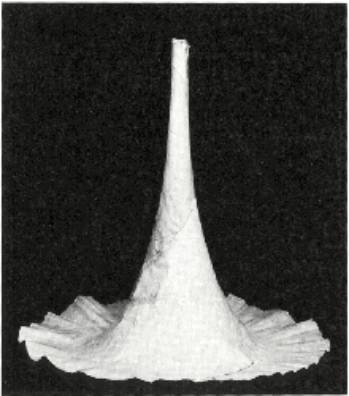}
        	\subcaption{}
        	\label{fig:1b}
\end{minipage}
~
\begin{minipage}[b]{0.32\textwidth}
	\includegraphics[width=\textwidth]{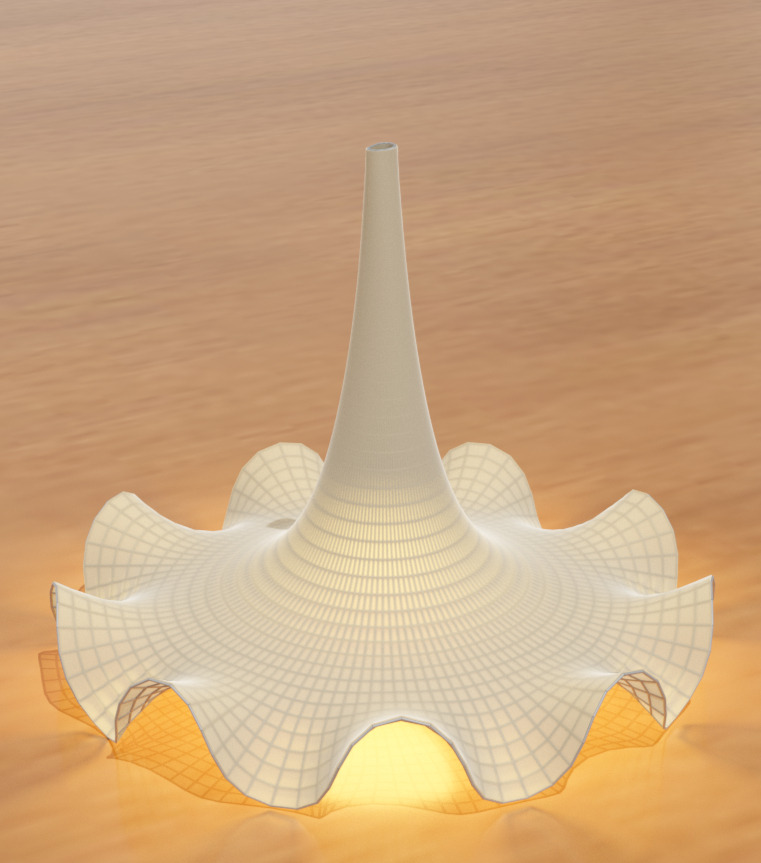}
        	\subcaption{}
        	\label{fig:1c}
\end{minipage}
\caption{
The pseudosphere: (a) as a region in the Poincar\'{e} disk model, (b) Beltrami's paper model (1869–1872, Dept.\ of Mathematics, Pavia), and (c) our mesh embedding.}
\label{fig:pseudosphere}
\end{figure}

\section*{Background and Related Work}

The desire to build accurate three-dimensional models is nearly as old as hyperbolic geometry itself. Shortly after its discovery, mathematicians sought concrete models to work with. Beltrami provided the first: the Beltrami-Klein disk and what we now call the Poincar\'{e} disk, which represent the infinite hyperbolic plane inside a bounded region. But like maps of the Earth, these models distort. Beltrami wanted the analog of a globe—a physical surface correctly realizing hyperbolic geometry. He succeeded in constructing paper models~\cite{EmmerAbate2020Beltrami}, including the pseudosphere (Figure~\ref{fig:pseudosphere}). There is an explicit formula for part of this surface, as a surface of revolution, but Beltrami's paper models extend further, ruffling and crinkling where equations give out.

Since Beltrami, others have approached this challenge using methods ranging from immersive visualizations and VR to physical and rendered models. Our work belongs to the latter category, and we briefly survey other works in this space. Henderson and Taimina's crocheted hyperbolic planes~\cite{henderson2001crocheting} offer a tactile, handmade alternative. Bulatov~\cite{bridges2011:479, bridges2013:167} and Fathauer~\cite{bridges2014:87} explore computational and paper folding approaches to hyperbolic tilings. Jackson and Williams~\cite{bridges2024:619} assemble Poincar\'{e} disk tessellations into three-dimensional models. Dunham~\cite{bridges2007:395} and Laigo et al.~\cite{bridges2008:311} analyze the geometry and symmetry of hyperbolic patterns. Physical models of related geometries appear in work by Demaine et al.~\cite{bridges1999:91}, Rousseau~\cite{bridges2003:127}, and Swart~\cite{bridges2015:151}. Kopczy\'{n}ski and Celi\'{n}ska~\cite{bridges2018:551} use energy minimization to embed surfaces for interactive VR exploration, a technique also used in HyperRogue~\cite{bridges2017:9}. Our underlying approach is similar, and we target rendered illustrations, 3D-printable models, and more general meshes with physical constraints for artistic applications. There is also prior work on embeddings of hyperbolic triangle tilings, often referred to as the Thurston model of the hyperbolic plane~\cite{bennett2010drawing}. Related constructions include Harriss’s physical polydron models~\cite{harriss_polydron}, hinged 3D-printed realizations by Irving and Segerman that appear in~\cite{segerman2016visualizing}, and paper models used in pedagogy and outreach, such as those of Henderson~\cite{hyperbolic_soccerball_iff} and Sottile~\cite{hyperbolic_soccerball_maa}.

Collectively, these works illustrate both the richness of and the technical subtleties involved in producing such embeddings. While local or structured constructions are possible, extending them to large regions while preserving intrinsic hyperbolic geometry is challenging, and no general method is known.

\section*{Creating Computer Models}

Our goal is to build models of regions in the hyperbolic plane (or surfaces built from hyperbolic geometry) as physical surfaces in $\mathbb{R}^3$. Explicit formulas for such embeddings are rare, so we need a simpler representation to compute with. Our approach is to discretize: we divide the surface into a \emph{mesh}, a network of vertices connected by edges, forming faces. Each edge has a length, each face an area, and they meet at vertices at measurable angles. This structure gives us a way to record the geometry of hyperbolic space in finite data.

In a model of the hyperbolic plane, we can compute these edge lengths exactly using the hyperbolic metric. An \emph{embedding} of the mesh into $\mathbb{R}^3$ assigns each vertex a position in space; the edges become line segments whose Euclidean lengths we can measure. The embedding is accurate when these measured lengths match the hyperbolic targets.  This gives us an easy way to check any proposed embedding. Just measure and compare. But it does not tell us how to construct one. The space of all possible vertex positions is vast, and most configurations are wrong.

\begin{figure}[h!tbp]
\centering
\begin{minipage}[b]{0.32\textwidth}
	\includegraphics[width=\textwidth]{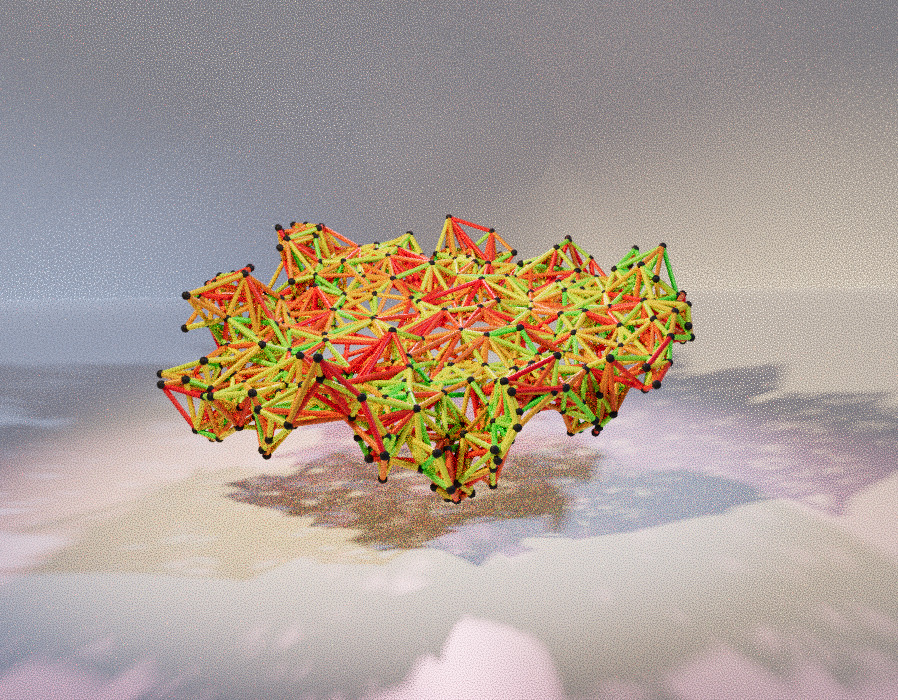}
\end{minipage}
~
\begin{minipage}[b]{0.32\textwidth}
	\includegraphics[width=\textwidth]{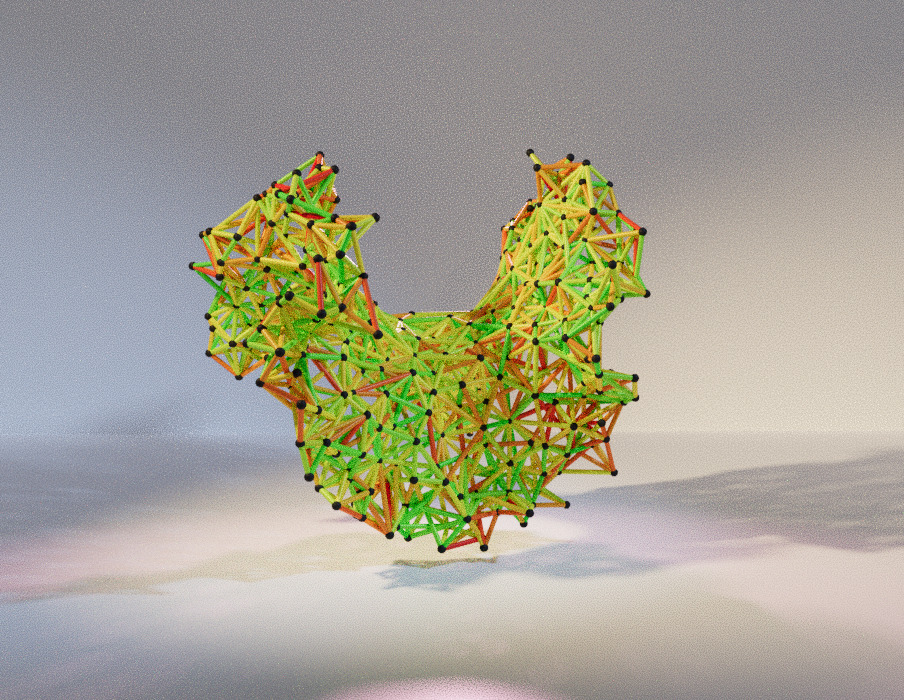}
\end{minipage}
~
\begin{minipage}[b]{0.32\textwidth}
	\includegraphics[width=\textwidth]{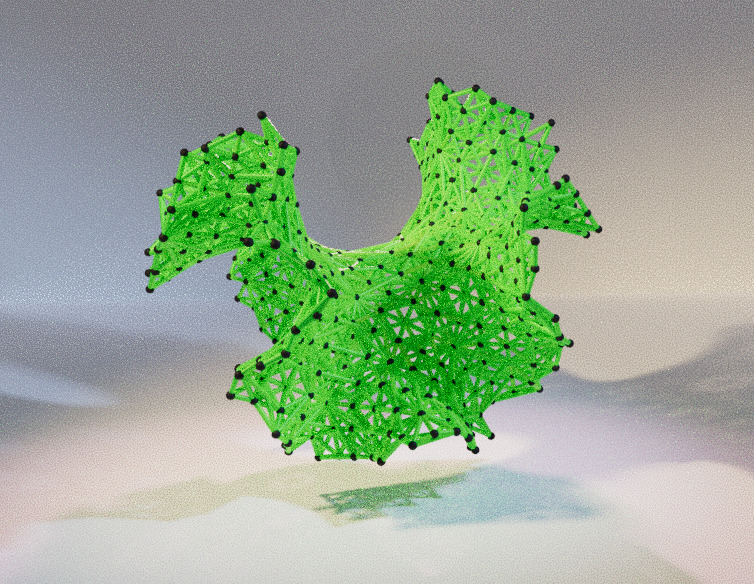}
\end{minipage}
\caption{A spring system encoding the geometry of a hyperbolic disk converging to an embedding from left to right.  Springs are colored red (poor) to green (good) depending on how well their length matches the corresponding hyperbolic length.}
\label{fig:embeddings}
\end{figure}

The solution is to search. We begin with an arbitrary embedding (typically a poor one, for instance a perturbed version of the disc model), define an energy measuring how far each edge deviates from its target length, and use gradient descent to iteratively improve. At each step, we nudge vertices in directions that reduce the total error. With enough iterations, the mesh settles into a configuration where the edge lengths are correct—an accurate model of the hyperbolic surface.
Figure~\ref{fig:embeddings} shows this process for a disk in the hyperbolic plane. Edges are colored by their deviation from the target length: red indicates significant error, green indicates agreement. The initial embedding is a crumpled mess; the final result is a smooth surface with the geometry of the hyperbolic plane built in. This places our method within the class of mass-spring systems, where the rest lengths are derived from an intrinsic (hyperbolic) metric.

We now make the method precise. Let $G=(V,E)$ be the edge graph (1-skeleton) of a mesh. An embedding of the mesh is a map
\begin{equation*}
x : V \to \mathbb{R}^3, \qquad i \mapsto x_i,
\end{equation*}
assigning a position $x_i \in \mathbb{R}^3$ to each vertex $i \in V$. The intrinsic geometry of the source surface determines a target length $\ell_{ij} > 0$ for each edge $(i,j)\in E$. Our goal is to find an embedding $x$ such that the Euclidean edge lengths $\|x_i - x_j\|$ match these prescribed values. We formulate this as an energy minimization problem. The primary term is a Hookean spring energy
\begin{equation*}
E_{\mathrm{spring}}(x)
=
\frac{1}{2}
\sum_{(i,j)\in E}
k_{ij}\,\bigl(\|x_i - x_j\| - \ell_{ij}\bigr)^2,
\end{equation*}
where $k_{ij} > 0$ are stiffness parameters. To discourage self-intersections, we include a Coulomb repulsion term
\begin{equation*}
E_{\mathrm{rep}}(x)
=
k_C \sum_{i < j} \frac{q_i q_j}{\|x_i - x_j\|},
\end{equation*}
where $q_i$ are vertex charges. In some scenes, we additionally incorporate optional gravity and positional constraints (e.g., fixed or pinned vertices). The total energy is
\begin{equation*}
E(x) = E_{\mathrm{spring}}(x) + E_{\mathrm{rep}}(x) \; (+\ \text{optional terms}),
\end{equation*}
and is minimized numerically. Starting from an initial embedding $x^{(0)}$ (e.g., a perturbed Poincar\'e disk layout, a cylinder, or a flat grid), we iteratively update vertex positions using gradient descent or damped second-order dynamics until convergence of the edge-length residuals. For example, the update rule for gradient descent is
\begin{equation*}
x^{(t+1)} = x^{(t)} - \eta \,\nabla E(x^{(t)}),
\end{equation*}
with step size $\eta > 0$. The framework is flexible with respect to mesh type, and we use both triangle and quad meshes. For quad meshes, we add diagonal (shear) and longer-range (bend) springs to stabilize the mesh embedding. For fixed initialization and parameters the procedure is deterministic, though some examples introduce small random perturbations or stochastic updates. The main user-controlled parameters are mesh resolution, stiffness coefficients, repulsion strength, time step, and iteration count.

\section*{Illustrating Hyperbolic Geometry}

Many models of hyperbolic space have been developed, from the aforementioned Beltrami-Klein and Poincar\'{e} disks to the upper half-plane, the hyperboloid and beyond. Each has its uses: the upper half-plane simplifies certain calculations, the hyperboloid lets us use linear algebra. However, they all drastically distort size. Embedded meshes complement these by making visible what models obscure.

\begin{figure}[h!tbp]
\centering
\begin{minipage}[b]{0.23\textwidth}
	\includegraphics[width=\textwidth]{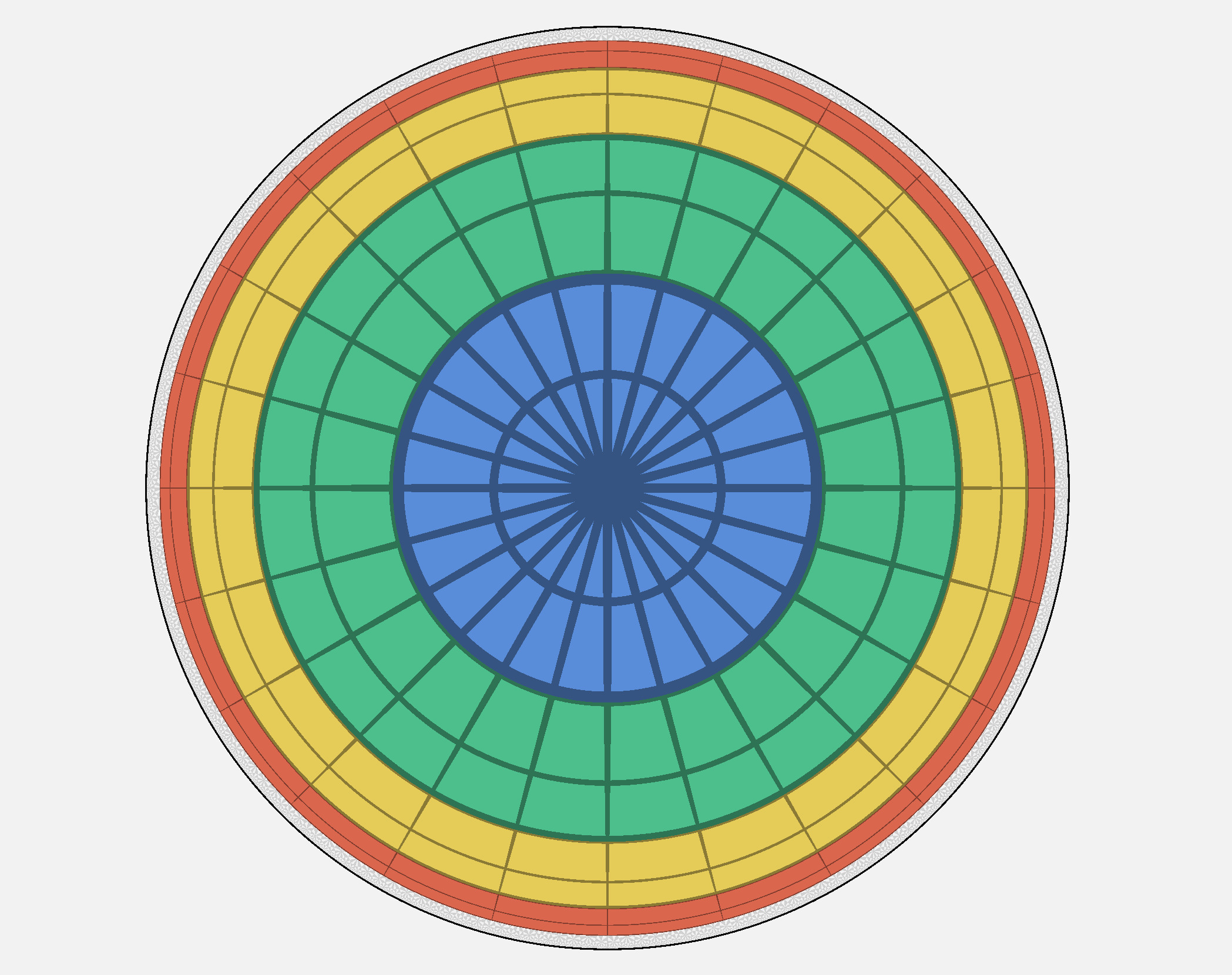}
    \subcaption{}
\end{minipage}
\begin{minipage}[b]{0.23\textwidth}
	\includegraphics[width=\textwidth]{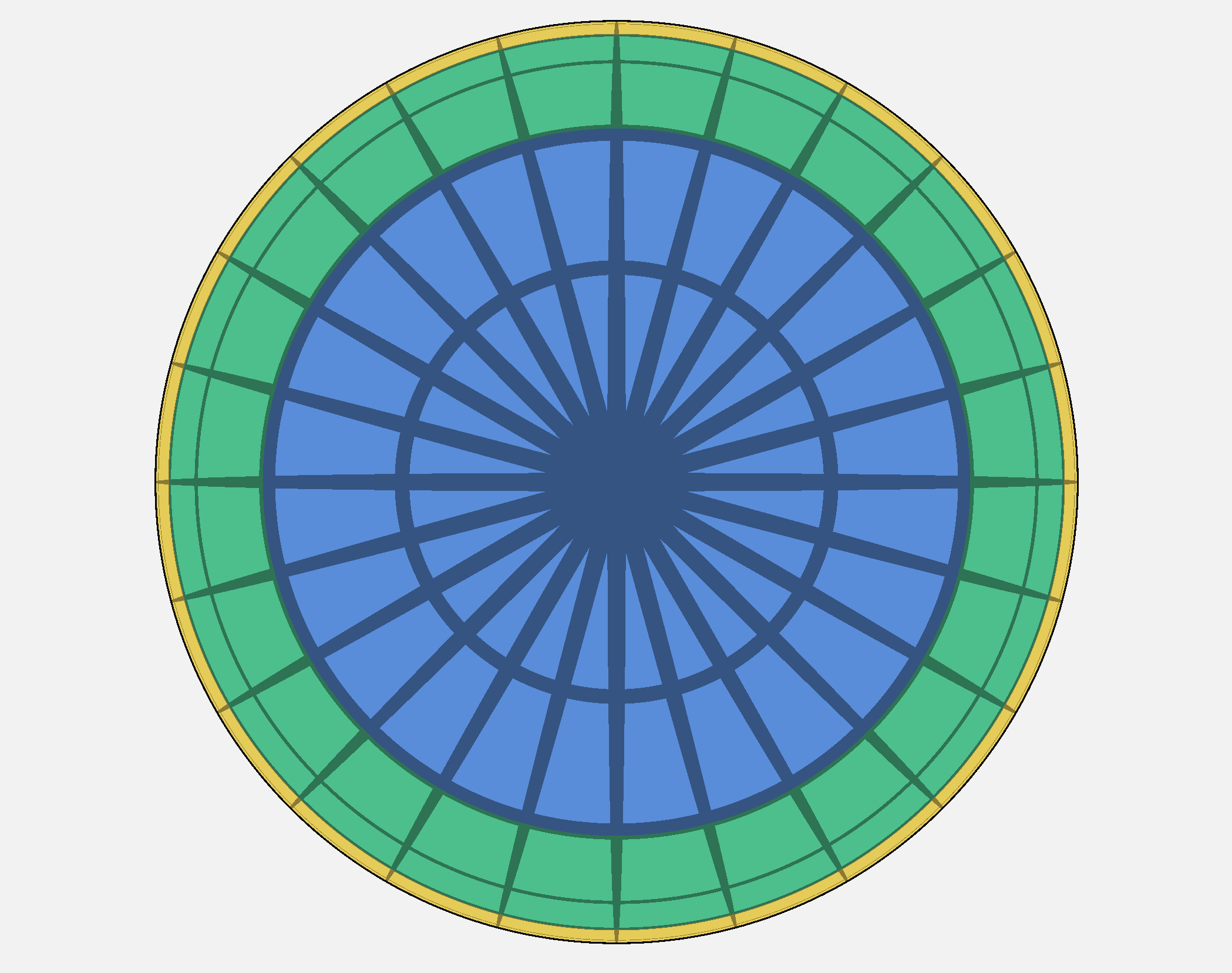}
    \subcaption{}
\end{minipage}
\begin{minipage}[b]{0.23\textwidth}
	\includegraphics[width=\textwidth]{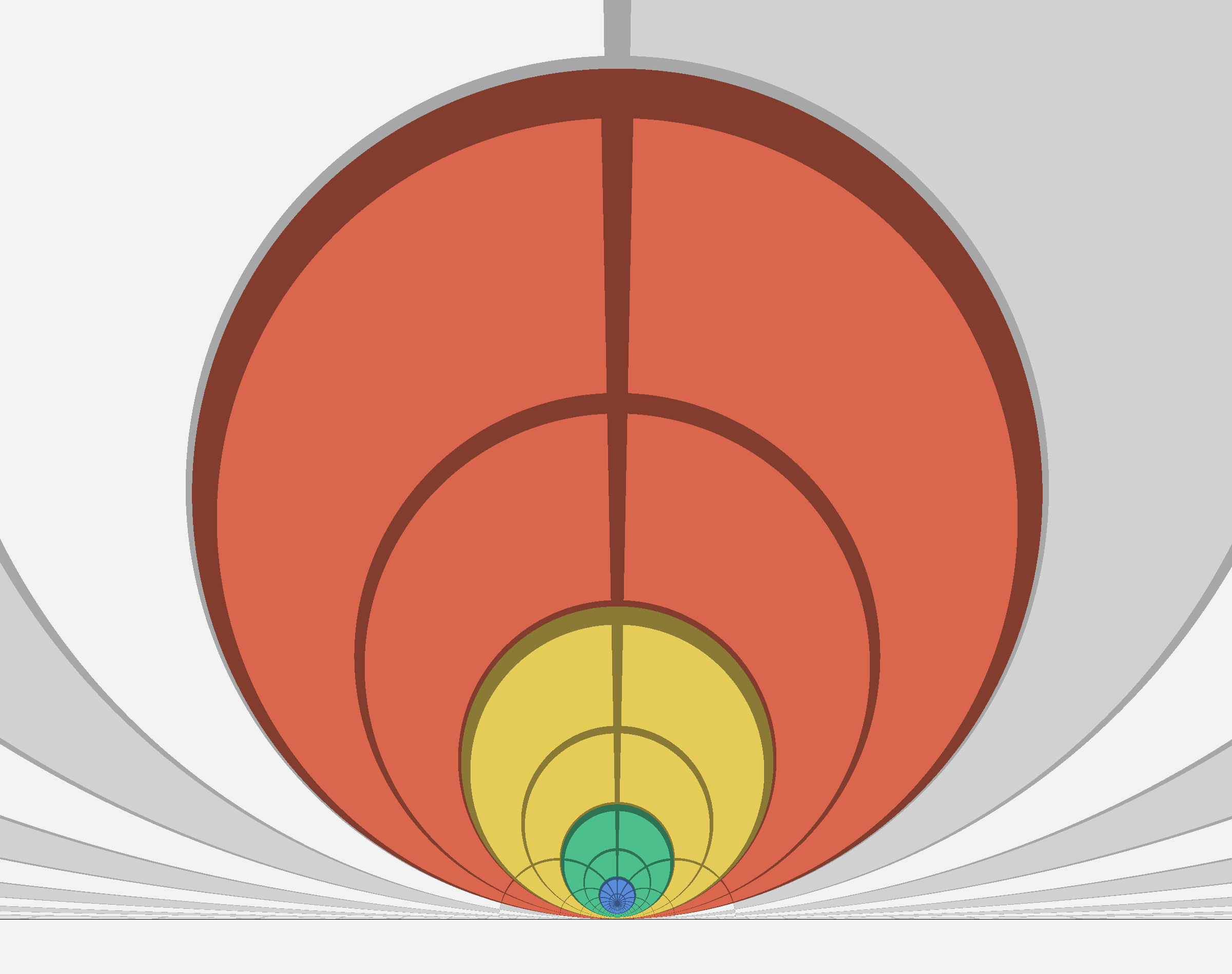}
    \subcaption{}
\end{minipage}
\begin{minipage}[b]{0.23\textwidth}
	\includegraphics[width=\textwidth]{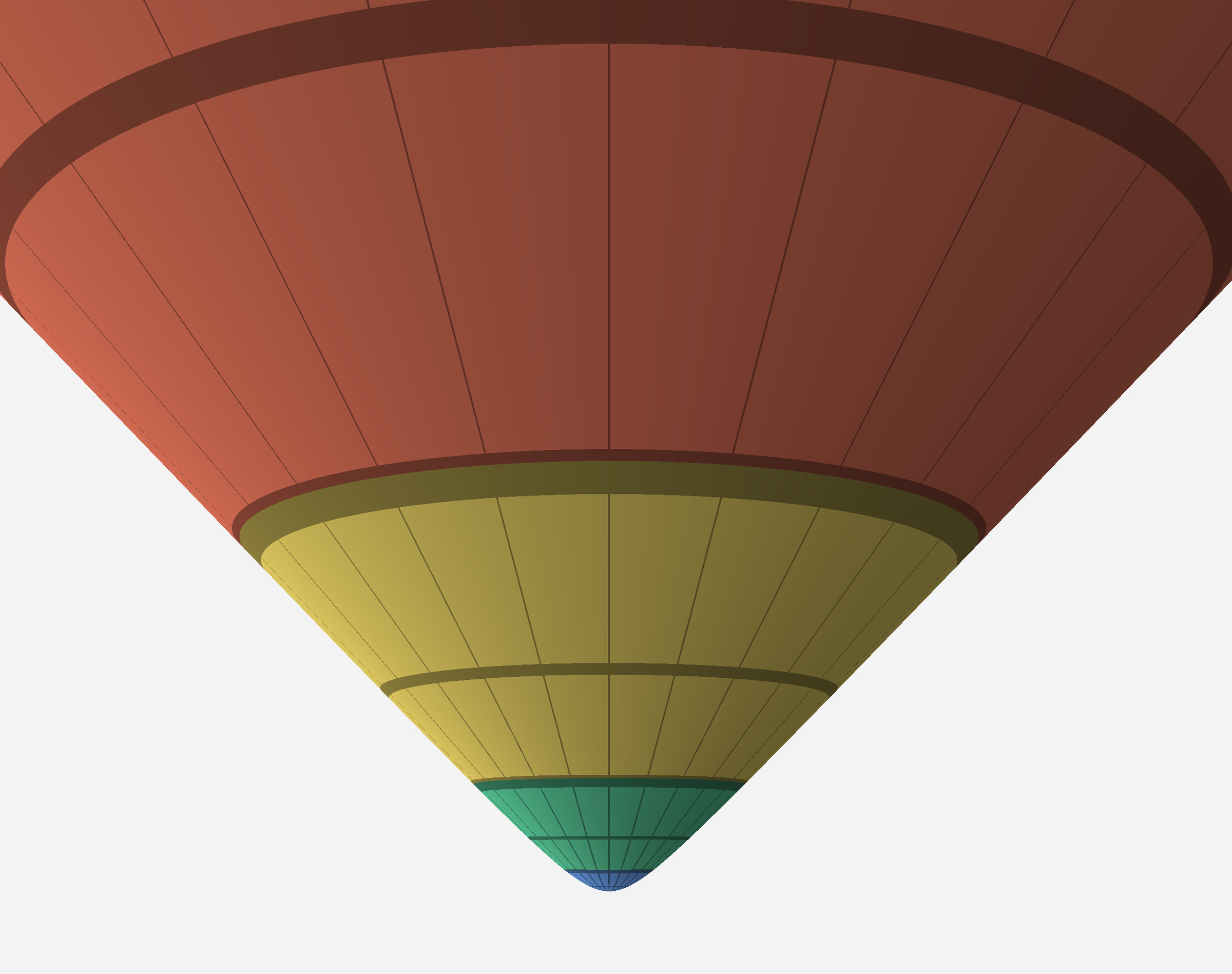}
    \subcaption{}
\end{minipage}

\caption{
Concentric hyperbolic disks of radius $1$ (blue), $2$ (green), $3$ (yellow), and $4$ (red), shown in (a) the Poincar\'{e} disk, (b) the Beltrami-Klein disk, (c) the upper half-plane,\\ and (d) the hyperboloid model.
}
\label{fig:disk-models}
\end{figure}

Figure~\ref{fig:disk-models} shows concentric disks of hyperbolic radius $1$, $2$, $3$, and $4$ in several models. In each, the disks appear to grow modestly—approaching the boundary in the disk models, climbing the hyperboloid, stretching asymmetrically in the half-plane. But hyperbolic area grows exponentially: a disk of radius $r$ has area $2\pi(\cosh(r)-1)\approx \pi e^{r}$ compared to $\pi r^2$ in Euclidean space.
The disk of radius $4$ is $48$ times as large as the disk of radius $1$! None of the models make this visible. But embedded in $\mathbb{R}^3$, explosive growth is clear.

The exponential growth of hyperbolic space causes real qualitative differences from Euclidean geometry (Figure \ref{fig:strip} in Appendix~\ref{app:extended-gallery-ii}). Consider a geodesic (the hyperbolic analog of a straight line) and the strip of points within some fixed distance of it. In Euclidean space this is just a rectangle: a strip of paper with the line down the center. In hyperbolic space, the curve at distance $d$ from the geodesic is $\cosh(d)$ times longer than the geodesic itself. The strip contains exponentially more material than its Euclidean counterpart. In the Poincar\'{e} disk this region appears as a lens or banana, seeming to contract toward the boundary; the exponential growth is invisible. Embedded in $\mathbb{R}^3$, the strip ruffles and buckles to accommodate its true area.

\begin{figure}[h!tbp]
\centering

\begin{minipage}[b]{0.23\textwidth}
	\includegraphics[width=\textwidth]{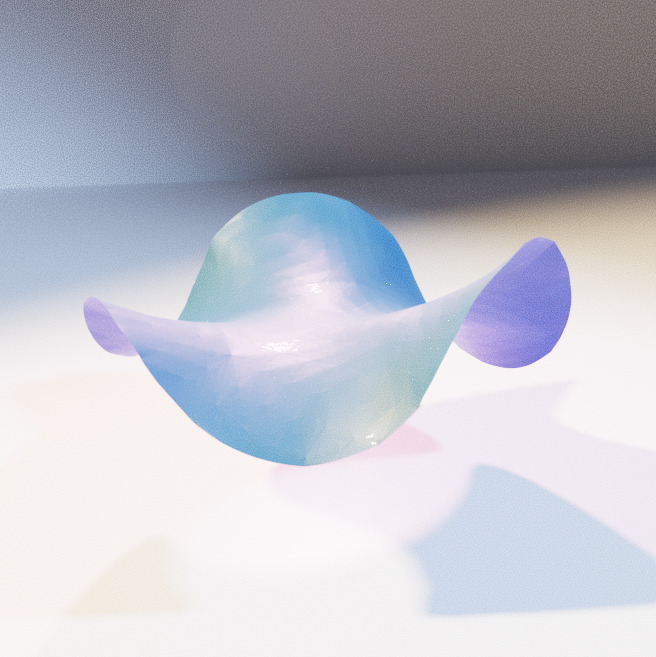}
\end{minipage}
\begin{minipage}[b]{0.23\textwidth}
	\includegraphics[width=\textwidth]{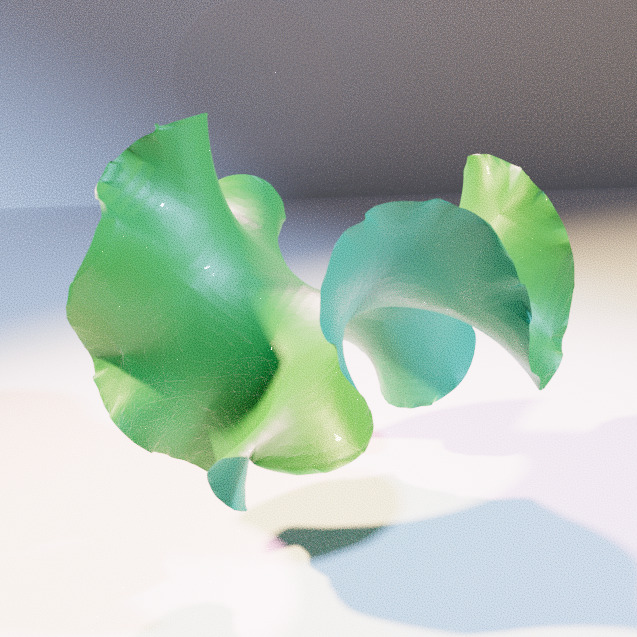}
\end{minipage}
\begin{minipage}[b]{0.23\textwidth}
	\includegraphics[width=\textwidth]{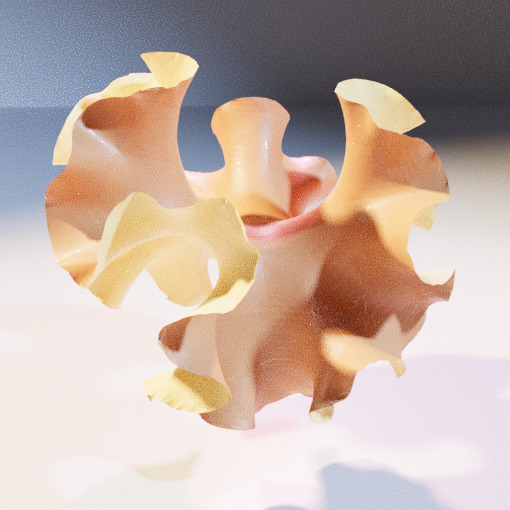}
\end{minipage}
\begin{minipage}[b]{0.23\textwidth}
	\includegraphics[width=\textwidth]{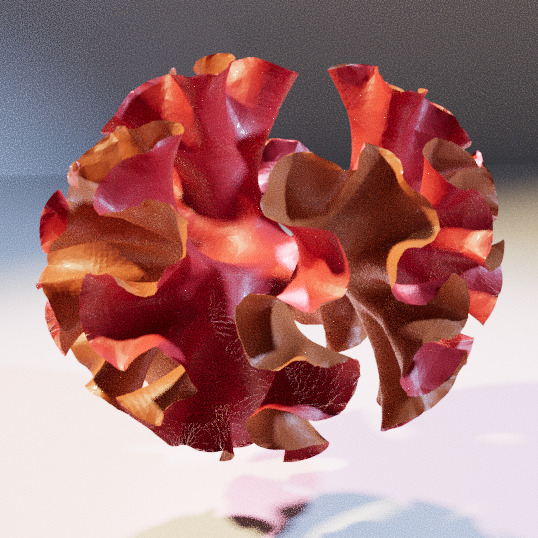}
\end{minipage}
\caption{Hyperbolic disks of radius
$1$ (blue), $2$ (green), $3$ (yellow), and $4$ (red),
embedded in Euclidean space. Compare with Figure~\ref{fig:disk-models}. }
\label{fig:disk-embeddings}
\end{figure}

Similarly stark is the behavior of parallel geodesics. In Euclidean space, parallel lines remain a fixed distance apart forever. In hyperbolic space, geodesics that start parallel are drawn apart exponentially quickly. The models show this divergence, but an embedding makes visceral just how fast space is growing between them.

\section*{Hyperbolic Geometry in Nature and Design}

Hyperbolic forms appear throughout the natural world. The ruffled edges of a lettuce leaf, the undulating fronds of kelp, the wrinkled surface of a coral, and the fluttering mantle of a cuttlefish are examples of such forms. All exhibit the characteristic buckling of surfaces with too much area to lie flat. This is no coincidence: biological growth often proceeds by adding material proportional to what is already there, and such processes lead inexorably to exponential expansion. The geometry discovered by Bolyai and Lobachevsky to probe the foundations of Euclid turns out to be encoded in organisms that have never heard of the parallel postulate.

\begin{figure}[h!tbp]
\centering

\begin{minipage}[b]{0.96\textwidth}
	\includegraphics[width=\textwidth]{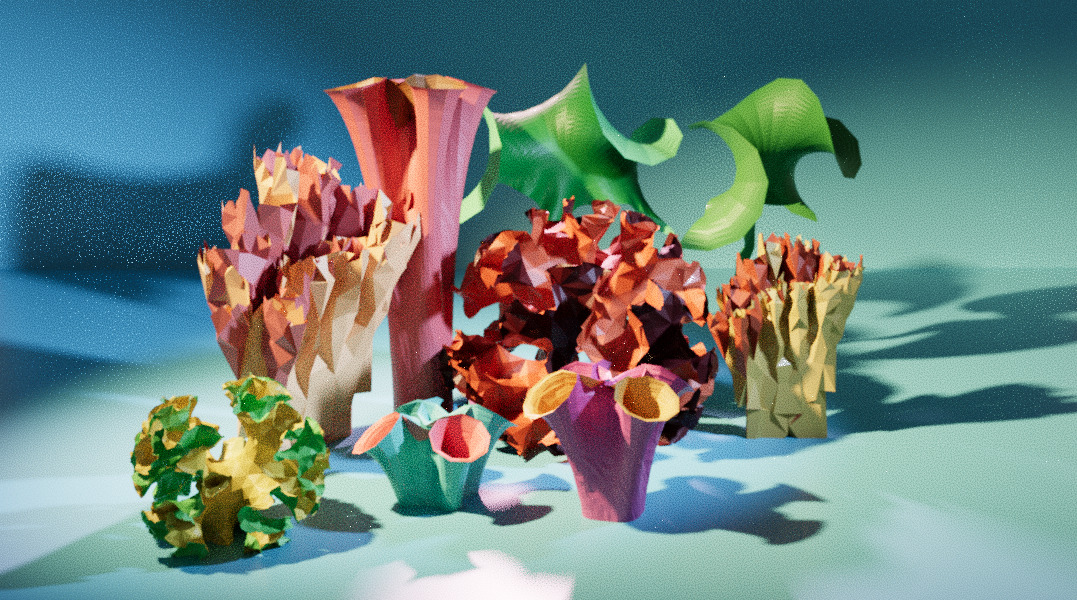}
        	\subcaption*{}
        	\label{fig:2d}
\end{minipage}
\caption{An artistic render using hyperbolic geometry to model natural forms in\\ a low-poly coral reef scene.}
\label{fig:nature}
\end{figure}

This connection is both beautiful and useful. Organic forms like these are difficult to model directly: no simple equation describes the edge of a lettuce leaf, and manually sculpting convincing ruffles is painstaking work. But an artist or designer seeking such forms can simply specify a hyperbolic mesh and embed it. The geometry does the work: energy minimization produces natural-looking ruffles automatically, without hand-tuning. The iterative optimization mirrors, in a poetic sense, the iterative process by which organisms grow: local adjustments accumulating into global form. Figure~\ref{fig:nature} shows examples from nature alongside a rendered scene built from embedded hyperbolic meshes. The coral is modeled not by artistic intuition but by geometry.

But geometry alone determines only shape, not orientation. In naturalistic art and design, the surface must be a participant in the world, responding to its environment: lettuce grows outward from its core, fabric hangs downward under gravity, and coral reaches for the light. The energy minimization framework accommodates this naturally. We can add terms beyond edge lengths: adding a gravitational potential
$mgh$ makes surfaces drape and settle. Adding pinning and collision constraints lets them hang from chosen points, rest on floors, or wrap around other objects. Figure~\ref{fig:dress} (Appendix~\ref{app:extended-gallery-ii}) shows a hyperbolic cylinder embedded without gravity (left) and with a gravitational term (center). The second hangs like the fabric of a ruffled dress. Some of the corals in Figure~\ref{fig:nature} were also generated this way, hanging a hyperbolic mesh under gravity and then rotating to upward-reaching forms.

\begin{figure}[h!tbp]
\centering
\includegraphics[width=0.96\textwidth]{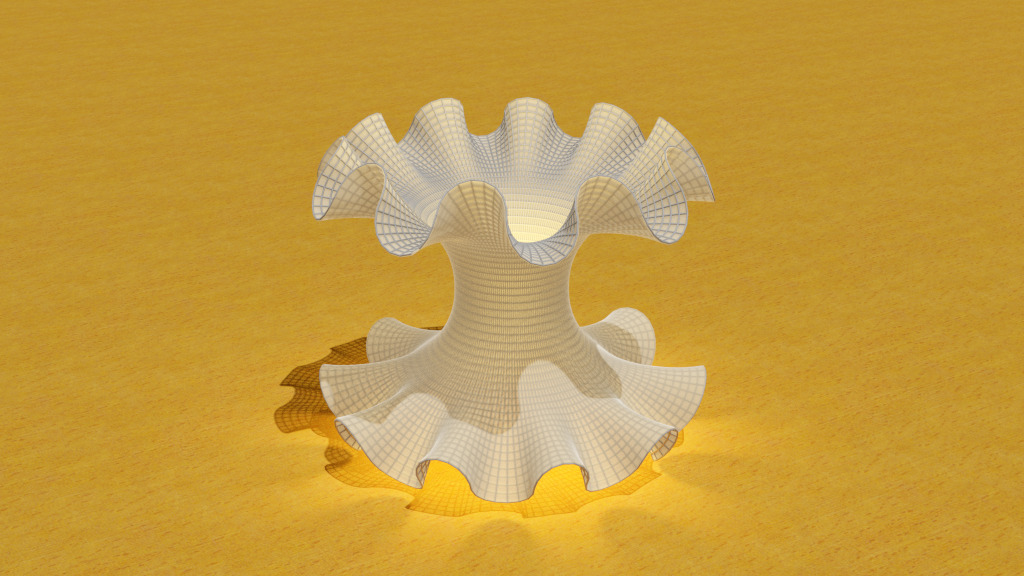}
\caption{
An embedded hyperbolic cylinder, rendered as a paper lantern with internal lighting.}
\label{fig:lantern}
\end{figure}

These techniques extend beyond modeling natural forms to design. Figure~\ref{fig:lantern} shows a hyperbolic cylinder rendered as a paper lantern. The intrinsic geometry alone produces its intricate form.

\section*{Utility in Public Engagement, Outreach, and Education}

\begin{figure}[h!tbp]
\centering
\hfill
\begin{minipage}[b]{0.48\textwidth}
	\includegraphics[width=\textwidth]{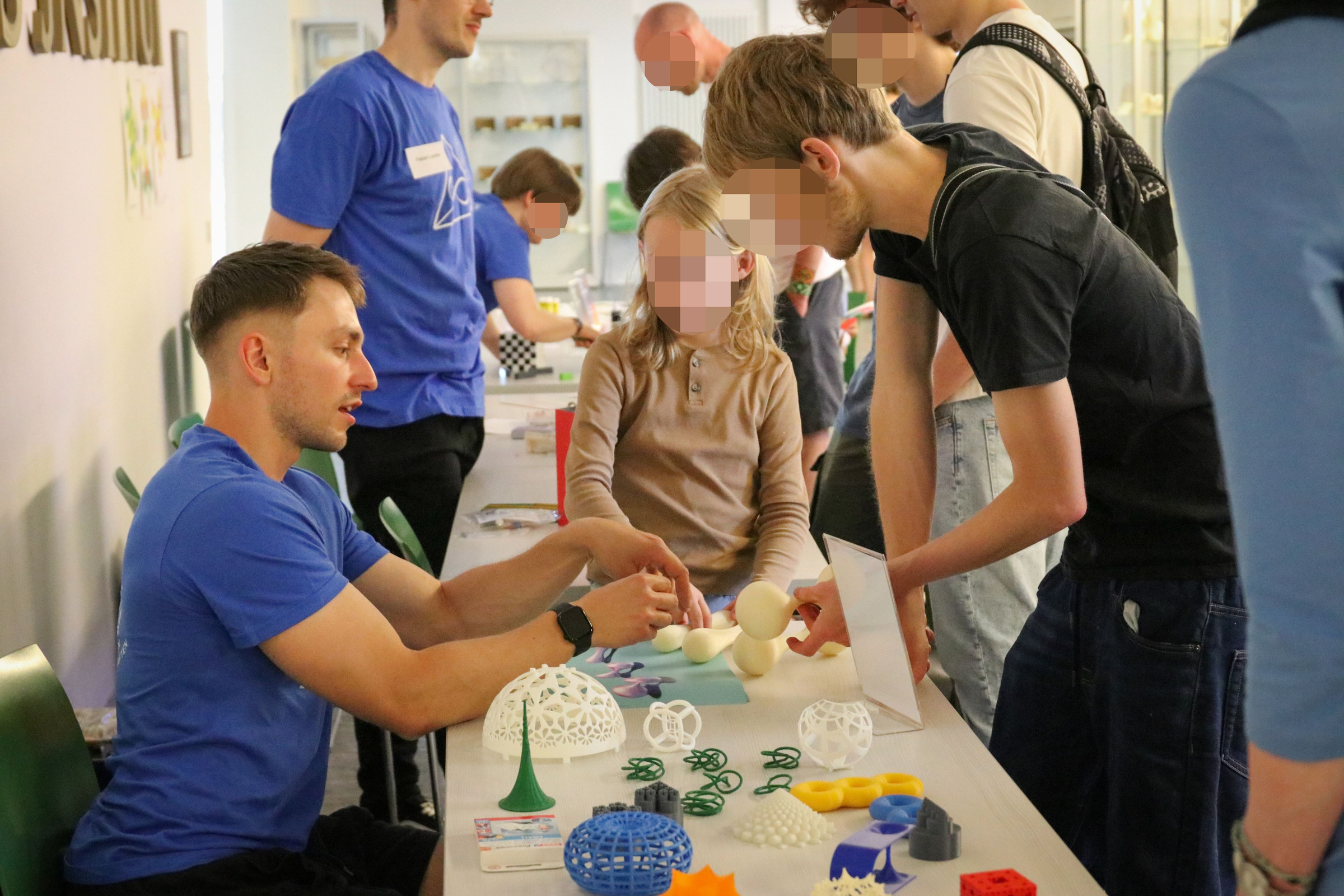}
\end{minipage}
\hfill
\begin{minipage}[b]{0.48\textwidth}
	\includegraphics[width=\textwidth]{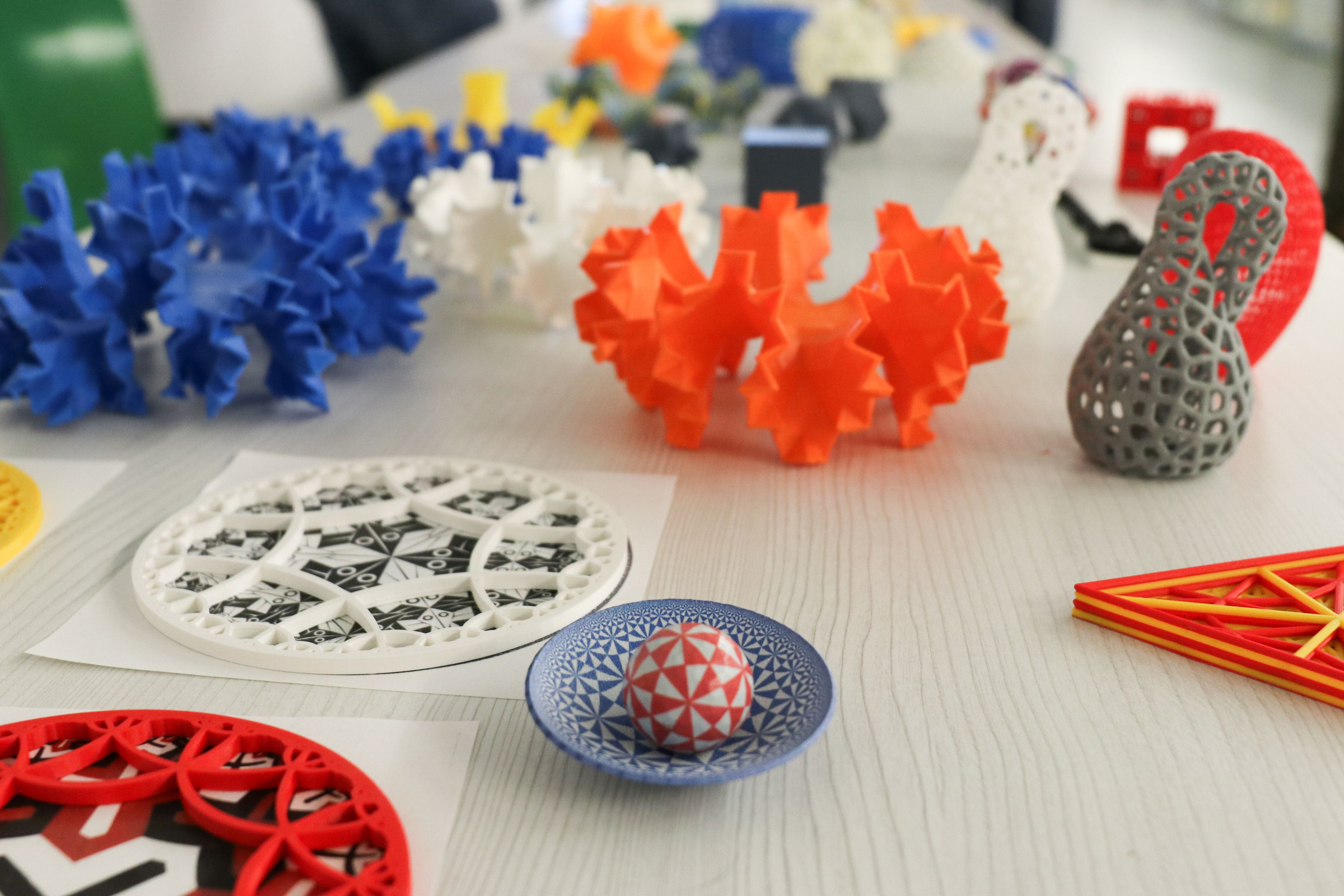}
\end{minipage}
\hfill
\caption{The first and second authors showcasing 3D-printed hyperbolic models at the Long Night of the Sciences, Leipzig University, 2025.}
\label{fig:outreach}
\end{figure}

We use our models, as renders or 3D prints, in research talks, public outreach events, school visits, art exhibits, and even as holiday postcards—among other settings. At talks, for example, displaying a render and rotating around a 3D model brings to life what flat models of hyperbolic geometry obscure: non-intersecting lines diverge, circles expand faster, and there is simply \emph{more space} than Euclidean intuition predicts. At public events (Figure~\ref{fig:outreach}), 3D prints and paper models invite hands-on exploration. Visitors can hold an embedding of a hyperbolic surface and see that its buckling is not decorative but necessary, the only way to fit the extra area into three-dimensional space. Producing some of these models has taken place as part of student and research projects, where building the embeddings requires engaging with both the underlying geometry and its computational implementation.

\section*{Conclusion}

We have shown how mesh embeddings can be used to visualize hyperbolic surfaces in Euclidean space, producing visualizations for illustration, art, and outreach. Future work could extend this work in several directions: producing illustrations of hyperbolic surfaces arising in research, embedding surfaces of higher genus, and exploring physical fabrication beyond 3D printing, such as fabric or paper templates.

\section*{Acknowledgements}

This work was carried out at the Mathematics Lab at the Max Planck Institute for Mathematics in the Sciences in Leipzig, Germany. The authors additionally acknowledge support of the Institut Henri Poincaré (UAR 839 CNRS-Sorbonne Université), and LabEx CARMIN (ANR-10-LABX-59-01). AI-assisted tools were used for minor editing tasks, including phrasing and clarity improvements. We sincerely thank the reviewers for their valuable comments and suggestions, which significantly improved the clarity and quality of this manuscript.

{\setlength{\baselineskip}{13pt}
\raggedright
\bibliographystyle{bridges}

\bibliography{references}

@inproceedings{bridges2024:619,
  author      = {Jackson, Thomas and Williams, Erin},
  title       = {Artfully Modeling Hyperbolic Planes through Tessellations},
  pages       = {619--626},
  booktitle   = {Bridges Conference Proceedings},
  year        = {Aug. 1--5, 2024},
  isbn        = {978-1-938664-49-6},
  issn        = {1099-6702},
  address     = {Richmond, Virginia, USA},

  url         = {http://archive.bridgesmathart.org/2024/bridges2024-619.html}
}

@inproceedings{bridges2014:87,
  author      = {Fathauer, Robert W.},
  title       = {Some Hyperbolic Fractal Tilings},
  pages       = {87--94},
  booktitle   = {Bridges Conference Proceedings},
  year        = {Aug. 14--19, 2014},
  isbn        = {978-1-938664-11-3},
  issn        = {1099-6702},
  address     = {Seoul, Korea},

  url         = {http://archive.bridgesmathart.org/2014/bridges2014-87.html}
}

@inproceedings{bridges2013:167,
  author      = {Bulatov, Vladimir},
  title       = {Bending Circle Limits},
  pages       = {167--174},
  booktitle   = {Bridges Conference Proceedings},
  year        = {Jul. 27--31, 2013},
  isbn        = {978-1-938664-06-9},
  issn        = {1099-6702},
  address     = {Enschede, the Netherlands},

  url         = {http://archive.bridgesmathart.org/2013/bridges2013-167.html}
}

@inproceedings{bridges2011:479,
  author      = {Bulatov, Vladimir},
  title       = {Bending Hyperbolic Kaleidoscopes},
  pages       = {479--482},
  booktitle   = {Bridges Conference Proceedings},
  year        = {Jul. 27--31, 2011},
  isbn        = {978-0-9846042-6-5},
  issn        = {1099-6702},
  address     = {Coimbra, Portugal},

  url         = {http://archive.bridgesmathart.org/2011/bridges2011-479.html}
}

@inproceedings{bridges2008:311,
  author      = {Laigo, Glenn R. and Puzon, Ia Kristine D. and Pe\~{n}as, Ma. Louise Antonette N. De Las},
  title       = {Coxeter Groups in Colored Tilings and Patterns},
  pages       = {311--318},
  booktitle   = {Bridges Conference Proceedings},
  year        = {Jul. 24--28, 2008},
  isbn        = {9780966520194},
  issn        = {1099-6702},
  address     = {Leeuwarden, the Netherlands},

  url         = {http://archive.bridgesmathart.org/2008/bridges2008-311.html}
}

@inproceedings{bridges2003:127,
  author      = {Rousseau, Irene},
  title       = {Geometric Mosaic Tiling on Hyperbolic Sculptures},
  pages       = {127--134},
  booktitle   = {Bridges Conference Proceedings},
  year        = {Jul. 23--25, 2003},
  isbn        = {84-930669-1-5},
  issn        = {1099-6702},
  address     = {Granada, Spain},

  url         = {http://archive.bridgesmathart.org/2003/bridges2003-127.html}
}

@inproceedings{bridges1999:91,
  author      = {Demaine, Erik D. and Demaine, Martin L. and Lubiw, Anna},
  title       = {Polyhedral Sculptures with Hyperbolic Paraboloids},
  pages       = {91--100},
  booktitle   = {Bridges Conference Proceedings},
  year        = {Jul. 30--Aug. 1, 1999},
  isbn        = {0-9665201-1-4},
  issn        = {1099-6702},
  address     = {Winfield, Kansas, USA},

  url         = {http://archive.bridgesmathart.org/1999/bridges1999-91.html}
}

@inproceedings{bridges2007:395,
  author      = {Dunham, Douglas},
  title       = {A `Circle Limit III' Calculation},
  pages       = {395--402},
  booktitle   = {Bridges Conference Proceedings},
  year        = {Jul. 24--27, 2007},
  isbn        = {0-9665201-8-1},
  issn        = {1099-6702},
  address     = {San Sebasti\'an, Spain},

  url         = {http://archive.bridgesmathart.org/2007/bridges2007-395.html}
}

@inproceedings{bridges2018:551,
  author      = {Kopczy\'{n}ski, Eryk and Celi\'{n}ska, Dorota},
  title       = {Virtual Crocheting of Euclidean Planes in a 3-Sphere},
  pages       = {551--554},
  booktitle   = {Bridges Conference Proceedings},
  year        = {Jul. 25--29, 2018},
  isbn        = {978-1-938664-27-4},
  issn        = {1099-6702},
  address     = {Stockholm, Sweden},

  url         = {http://archive.bridgesmathart.org/2018/bridges2018-551.html}
}

@inproceedings{bridges2015:151,
  author      = {Swart, David},
  title       = {Soccer Ball Symmetry},
  pages       = {151--158},
  booktitle   = {Bridges Conference Proceedings},
  year        = {Jul. 29--Aug. 1, 2015},
  isbn        = {978-1-938664-15-1},
  issn        = {1099-6702},
  address     = {Baltimore, Maryland, USA},

  url         = {http://archive.bridgesmathart.org/2015/bridges2015-151.html}
}

@incollection{EmmerAbate2020Beltrami,
  author    = {Andreatta, Marco},
  title     = {The Rise of Abstractionism: Art and Mathematics},
  booktitle = {Imagine Math 8},
  editor    = {Emmer, Michele and Abate, Marco},
  publisher = {Springer},
  year      = {2020}
}

@article{henderson2001crocheting,
  title={Crocheting the Hyperbolic Plane},
  author={Henderson, David W and Taimina, Daina},
  journal={Mathematical Intelligencer},
  volume={23},
  number={2},
  pages={17--27},
  year={2001},
  publisher={Springer-Verlag 175 FIFTH AVE, NEW YORK, NY 10010 USA}
}

@inproceedings{bridges2017:9,
  author      = {Kopczy\'{n}ski, Eryk and Celi\'{n}ska, Dorota and \v{C}trn\'{a}ct, Marek},
  title       = {HyperRogue: Playing with Hyperbolic Geometry},
  pages       = {9--16},
  booktitle   = {Bridges Conference Proceedings},
  year        = {Jul. 27--31, 2017},
  isbn        = {978-1-938664-22-9},
  issn        = {1099-6702},
  address     = {Waterloo, Ontario, Canada},

  url         = {http://archive.bridgesmathart.org/2017/bridges2017-9.html}
}

@article{bennett2010drawing,
  title={Drawing a Triangle on the Thurston Model of Hyperbolic Space},
  author={Bennett, Curtis D and Mellor, Blake and Shanahan, Patrick D},
  journal={Mathematics Magazine},
  volume={83},
  number={2},
  pages={83--99},
  year={2010},
  publisher={Taylor \& Francis}
}

@book{segerman2016visualizing,
  title={Visualizing Mathematics with 3D Printing},
  author={Segerman, Henry},
  year={2016},
  publisher={JHU Press}
}

@webpage{harriss_polydron,
  author = {Edmund Harriss},
  title = {\emph{Maxwell's Demon}},
  howpublished = {Hyperbolic Polydron},
  year = {2009},
  url = {https://maxwelldemon.com/2009/04/13/unscheduled-post-hyperbolic-polydron/}
}

@webpage{hyperbolic_soccerball_iff,
  author = {Keith Henderson},
  howpublished = {How to Build Your Own Hyperbolic Soccer Ball Model},
  year = {n.d.},
  url = {https://www.theiff.org/images/IFF_HypSoccerBall.pdf}
}

@webpage{hyperbolic_soccerball_maa,
  author = {Frank Sottile},
  title={\emph{Frank Sottile's Personal Website}},
  howpublished = {Make a Hyperbolic Football!},
  year = {2019},
  url = {https://franksottile.github.io/research/stories/hyperbolic_football/index.html}
}
}

\clearpage
\appendix

\section{Extended Gallery I}
\label{app:extended-gallery-1}

\begin{figure}[h!tbp]
\centering
\begin{minipage}[b]{0.23\textwidth}
	\includegraphics[width=\textwidth]{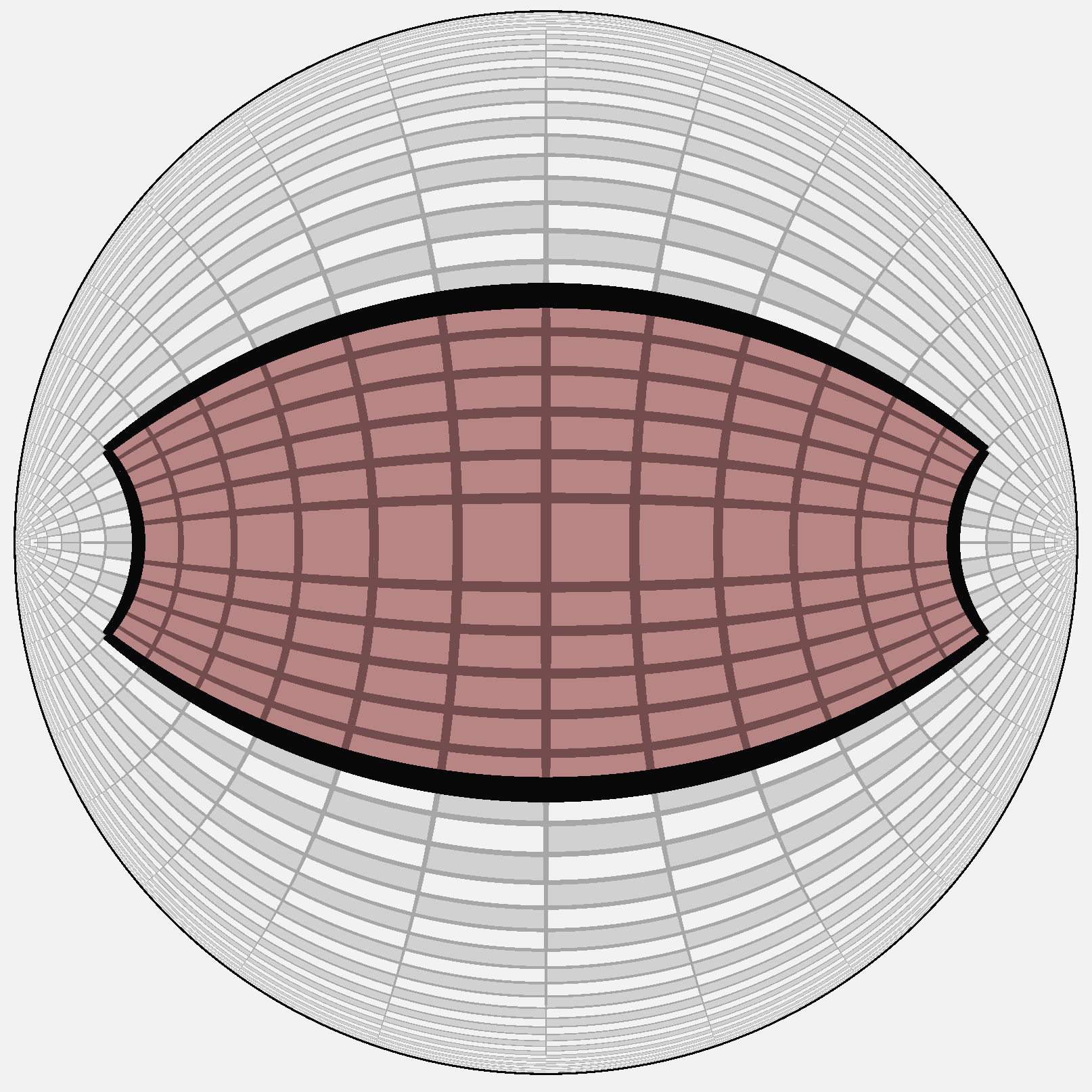}
\end{minipage}
\begin{minipage}[b]{0.23\textwidth}
	\includegraphics[width=\textwidth]{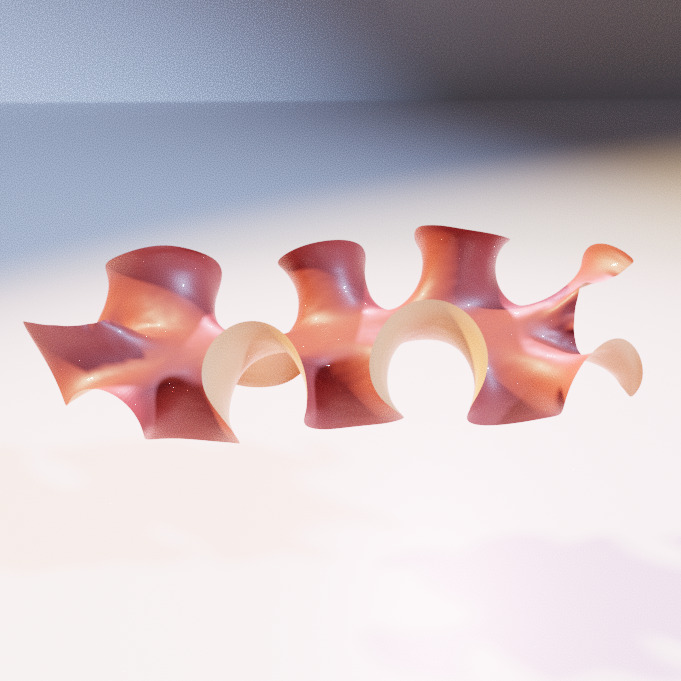}
\end{minipage}
\quad
\begin{minipage}[b]{0.23\textwidth}
	\includegraphics[width=\textwidth]{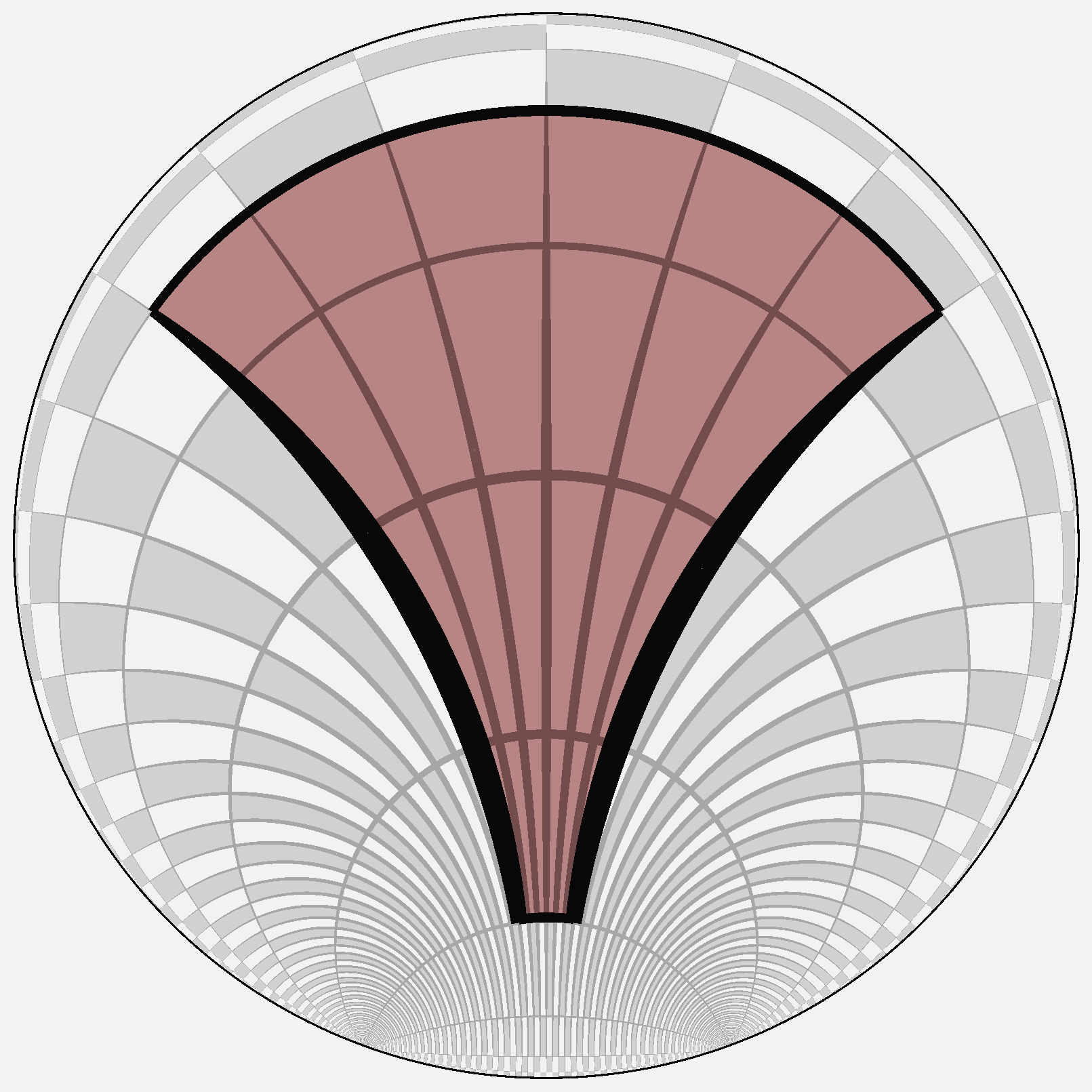}
\end{minipage}
\begin{minipage}[b]{0.23\textwidth}
	\includegraphics[width=\textwidth]{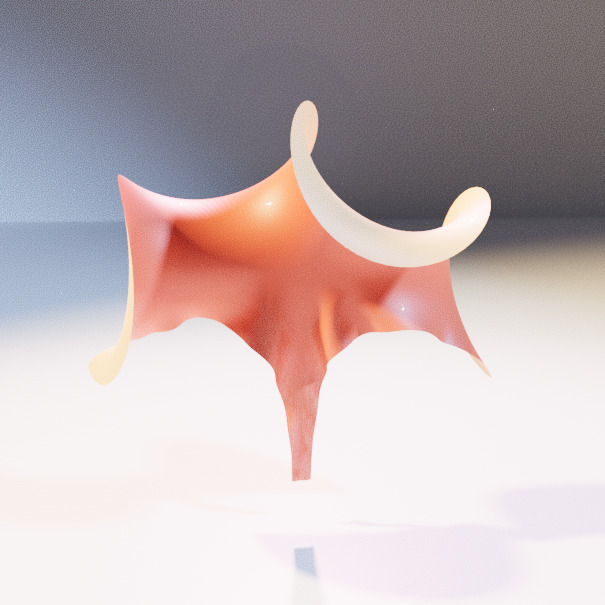}
\end{minipage}
\caption{A strip of points within a fixed distance to a geodesic in the hyperbolic plane, left.  A spray of geodesics all initially parallel, right.}
\label{fig:strip}
\end{figure}

\begin{figure}[h!tbp]
\centering
\begin{minipage}[b]{0.32\textwidth}
	\includegraphics[width=\textwidth]{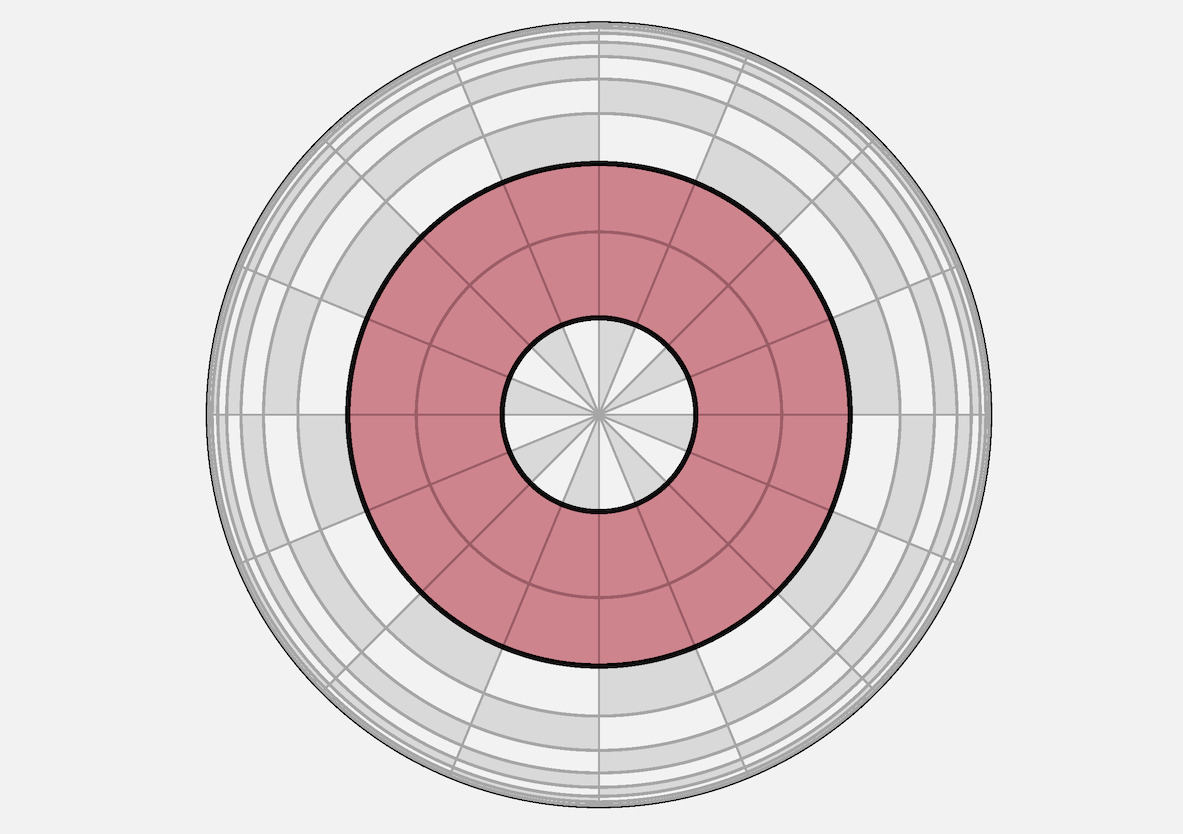}
    \subcaption{}
\end{minipage}
~
\begin{minipage}[b]{0.32\textwidth}
	\includegraphics[width=\textwidth]{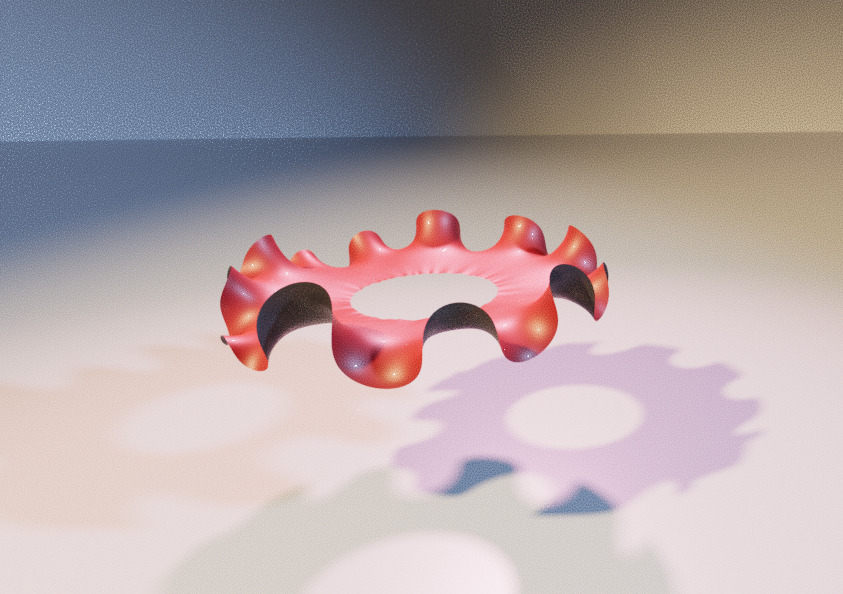}
    \subcaption{}
\end{minipage}
~
\begin{minipage}[b]{0.32\textwidth}
	\includegraphics[width=\textwidth]{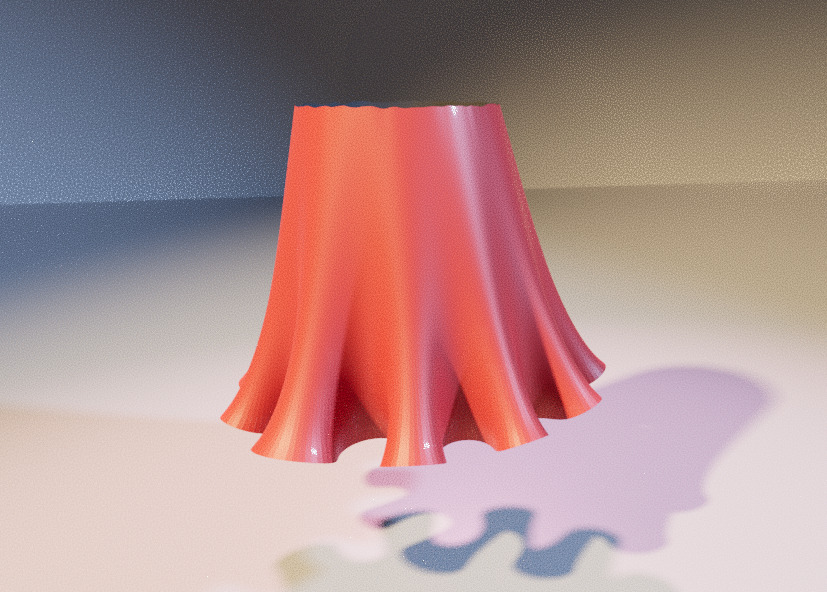}
    \subcaption{}
\end{minipage}
\caption{A ruffled dress modeled with hyperbolic geometry: (a) the initial region in the Poincar\'{e} disk, (b) the embedded mesh without gravity, (c) the embedded mesh, relaxed under gravity.}
\label{fig:dress}
\end{figure}

The construction process for the models in Figure~\ref{fig:nature} is illustrated in Figures~\ref{fig:strip} and~\ref{fig:dress}. In each case, we begin with a region of the hyperbolic plane (such as a strip or disk), discretize it as a mesh with edge lengths determined by the intrinsic metric, and compute an embedding in $\mathbb{R}^3$ via the spring-based optimization described above. Throughout, we use well-known analytic expressions for constructing circles, geodesics, and horocycles, and for computing distances in the hyperbolic plane. The resulting surfaces are then rendered with simple geometric or artistic interpretations, producing the forms shown in Figure~\ref{fig:nature}.

\section{Extended Gallery II}
\label{app:extended-gallery-ii}

\begin{figure}[h!tbp]
\centering
\includegraphics[width=\textwidth]{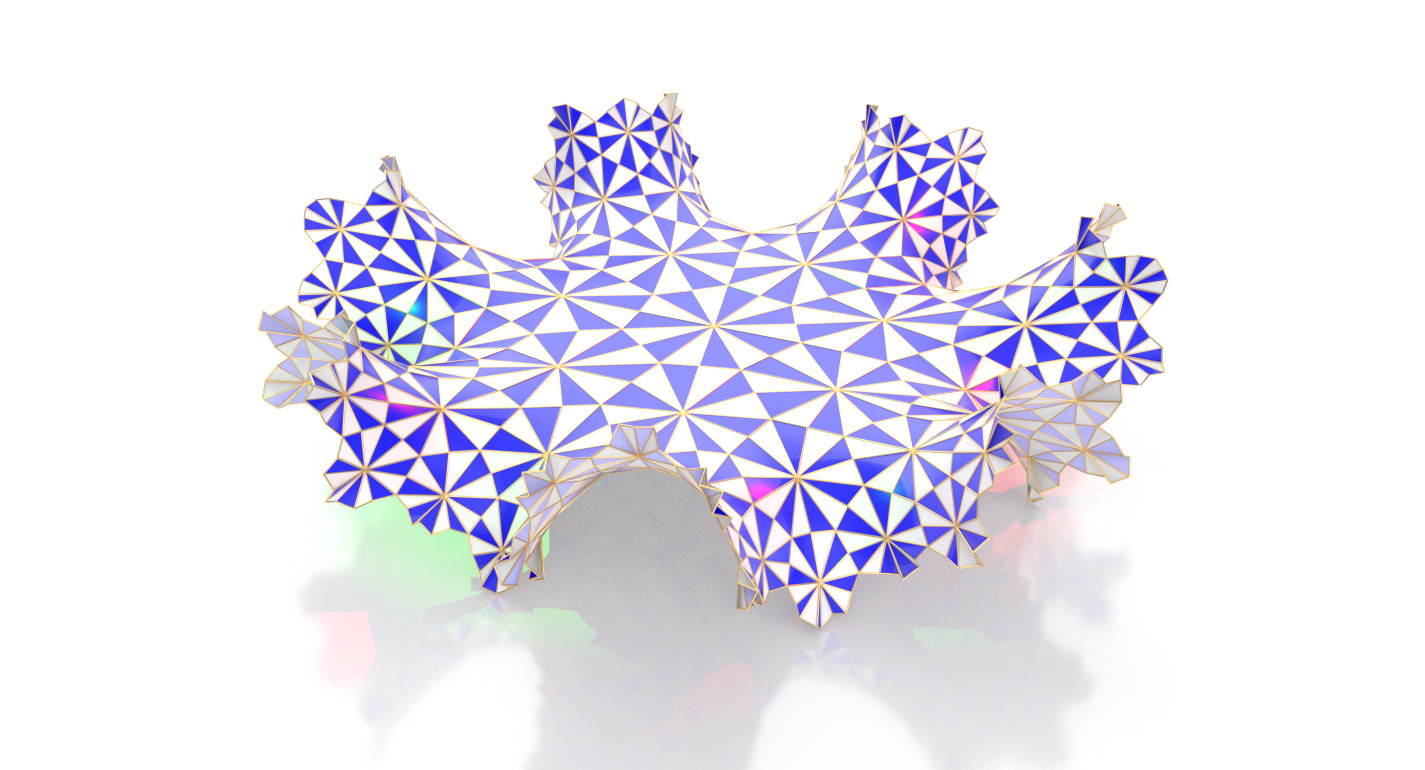}
\caption{
An embedding of a portion of a $(2,3,7)$ triangular tiling of the hyperbolic disk.}
\label{fig:disk-2-3-7}
\end{figure}

\begin{figure}[h!tbp]
\centering
\includegraphics[width=\textwidth]{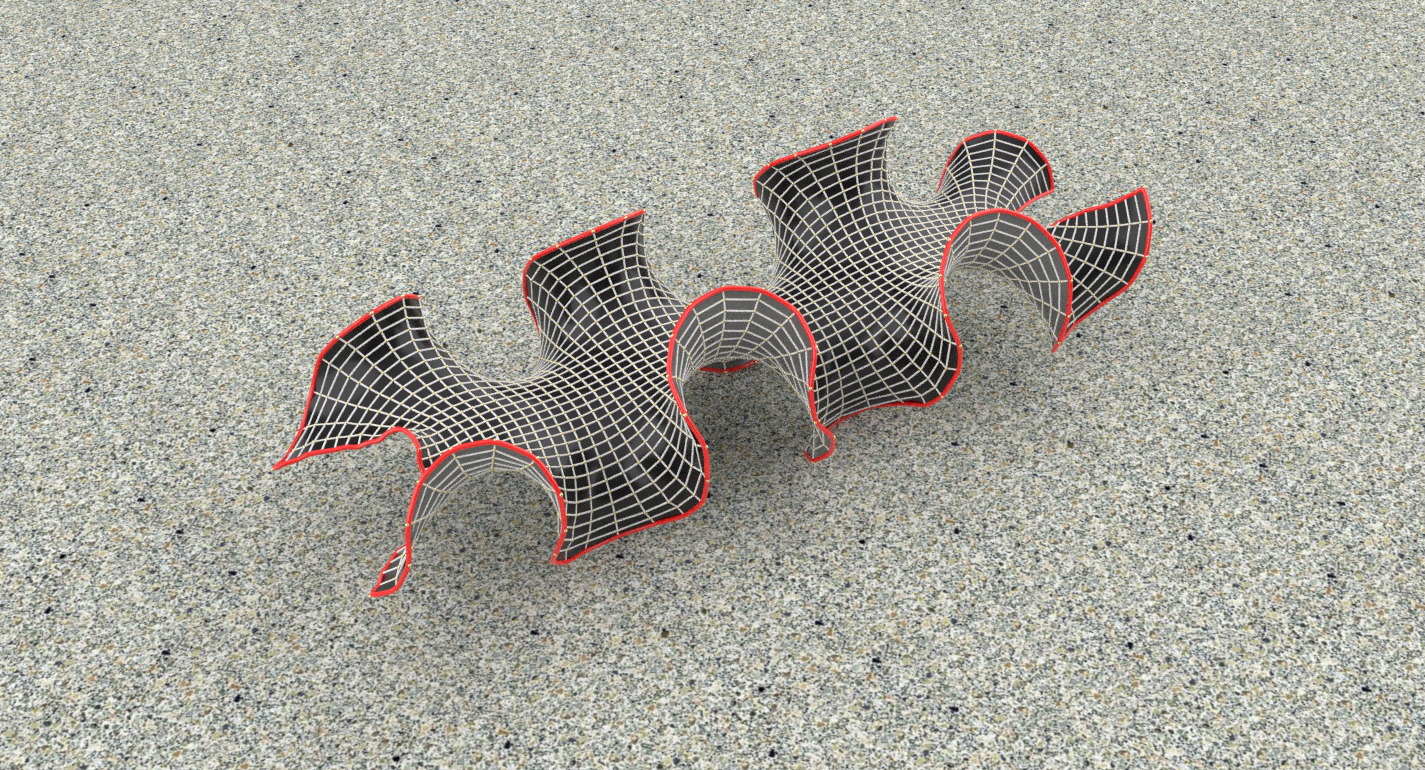}
\caption{
An embedded hyperbolic strip,  rendered as a sea slug against a terrazzo floor.}
\label{fig:strip-seaslug}
\end{figure}

\begin{figure}[h!tbp]
\centering
\includegraphics[width=\textwidth]{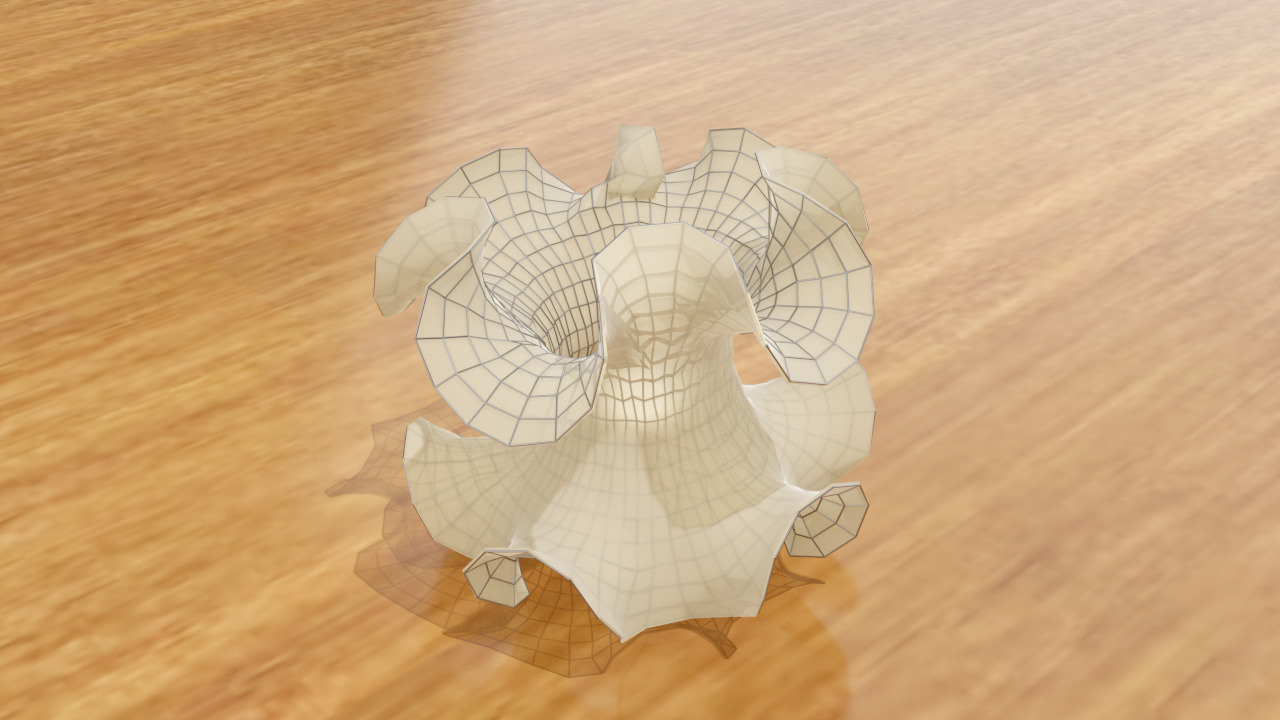}
\caption{
An embedded low-poly hyperbolic cylinder, rendered as a paper lantern with internal lighting. The optimization was terminated before full convergence.}
\label{fig:lantern-lowpoly}
\end{figure}

In this appendix, we include a selection of early experimental renderings not shown in the main text. These include a $(2,3,7)$ triangular tiling of the hyperbolic disk (Figure~\ref{fig:disk-2-3-7}), a hyperbolic strip interpreted as a sea slug (Figure~\ref{fig:strip-seaslug}), and a hyperbolic cylinder interpreted as a Japanese paper lantern (Figures~\ref{fig:lantern-lowpoly} and~\ref{fig:lantern}). These examples reflect the authors' early exploratory and creative process in developing the models for this work.

\end{document}